\documentclass[a4paper,11pt]{amsart}

\usepackage{nag}
\usepackage[utf8]{inputenc}
\usepackage[T1]{fontenc}
\usepackage{amsmath,amsthm,latexsym,amssymb,enumerate,accents}
\usepackage{color}
\usepackage[english]{babel}
\usepackage[colorlinks=false, pdfborder={0 0 0}]{hyperref}
\usepackage[left=3cm,right=3cm,top=2cm,bottom=3cm]{geometry}
\usepackage{booktabs}
\usepackage{graphicx}
\hypersetup{pageanchor=false}
\hypersetup{colorlinks,linkcolor={blue},citecolor={cyan},urlcolor={magenta}}  

\usepackage{subcaption}
\usepackage{url}
\usepackage{tikz}

\newtheorem{Thm}{Theorem}[section]

\newtheorem{Defi}[Thm]{Definition}
\newtheorem{Rma}[Thm]{Remark}

\definecolor{darkgreen}{RGB}{0,100,0}
\definecolor{darkred}{RGB}{180,0,0}
\definecolor{darkblue}{RGB}{0,0,150}
\definecolor{mustard}{RGB}{180,120,0}
\definecolor{zurzuvi}{RGB}{200,0,200}

\DeclareMathOperator{\kernel}{Ker}
\DeclareMathOperator{\im}{Im}

\DeclareMathOperator{\R}{\mathbb{R}}
\DeclareMathOperator{\DFlag}{dFl}
\DeclareMathOperator{\F}{\mathbb{F}}
\newcommand{\hadgesh}[1]{\emph{\textcolor{blue}{#1}}}

\newcommand{\modifprev}[1]{\textcolor{black}{#1}}

\title{\bfseries Computing persistent homology of directed flag complexes}

\begin{document}
\begin{abstract}
We present a new computing package \textsc{Flagser}, designed to construct the directed  flag complex of a finite directed graph, and compute  persistent homology for flexibly defined filtrations on the graph and the resulting complex.  The persistent homology computation part  of \textsc{Flagser} is based on the program \textsc{Ripser} \cite{Bauer}, but is optimised specifically for large computations. The construction of the directed flag complex is done in a way that allows easy parallelisation by arbitrarily many cores.  \textsc{Flagser} also has the option of working with undirected graphs. For homology computations \textsc{Flagser} has an  Approximate option, which shortens  compute time with remarkable accuracy. We demonstrate the power of \textsc{Flagser} by applying it to the construction of the directed flag complex  of digital reconstructions of brain microcircuitry by the Blue Brain Project and several other examples. In some instances we perform computation of homology.  For a more complete performance analysis, we also apply \textsc{Flagser} to some other data collections. In all cases the hardware used in the computation, the use of memory and the compute time are recorded.
\end{abstract}

\author{D. L\"utgehetmann}
\author{D. Govc}
\author{J.P. Smith}
\author{R. Levi}
\maketitle

%%%%%%%%%%%%%%%%%%%%%%%%%%%%%%%%%%%%%%%%%
%%%%%%%%%%%%%%%%%%%%%%%%%%%%%%%%%%%%%%%%%
%%%%%%%%%%%%%%%%%%%%%%%%%%%%%%%%%%%%%%%%%
%%%%%%%%%% INTODUCTION      %%%%%%%%%%%%%%%%%%%%%
%%%%%%%%%%%%%%%%%%%%%%%%%%%%%%%%%%%%%%%%%
%%%%%%%%%%%%%%%%%%%%%%%%%%%%%%%%%%%%%%%%%
%%%%%%%%%%%%%%%%%%%%%%%%%%%%%%%%%%%%%%%%%

\section{Introduction}
Over the past decade, homology and persistent homology --- particularly in the form of persistent homology of Vietoris-Rips complexes of point clouds --- entered the field of data analysis.
Therefore, its efficient computation is of great interest, which led to several computational tools like \textsc{Gudhi} \cite{Gudhi}, \textsc{PHAT} \cite{PHAT} and other packages.

In  an ongoing collaboration with the Blue Brain Project (Link \ref{BBP}) and the Laboratory for Topology and Neuroscience (Link \ref{UPHESS}) we study certain ordered simplicial complexes (Definition \ref{Def-Directed-Flag-Complex}) arising from directed graphs that model  brain microcircuitry reconstructions created by the Blue Brain  Project team. Having access to these  reconstructions gives us a unique opportunity to apply topological tools to the study of  brain circuitry models in a scale and level of biological accuracy that was never possible before. Thus the need arose for a software package that can generate the complexes in question out of the adjacency matrix of the circuits under investigation  and compute homology and persistent homology of such complexes with respect to a variety of filtrations. 

The size of the networks we considered initially was at the order of magnitude of 30 thousand vertices and 8 million directed edges (reconstructed neocortical column of a 14 days old rat (Link \ref{BBP-Rat})). Software was designed as part of the work  that lead to the publication \cite{Reimann17}, and is freely available (Link \ref{Neurotop}). However, the magnitude of the resulting complexes  - typically 7 dimensional with about 70 million 2- and 3-simplices - limited our ability to conduct certain deeper investigations, as even the construction of the complexes required a considerable number of compute hours on a large HPC node, and  homology computations with the full complex were completely out of reach. More recent reconstructions by the BBP team produced much larger and more complex models (See for instance the model of the full mouse cortex: Link \ref{BBP-Mouse}). In addition, attempting to understand neuroscientific phenomena, such as neuronal activity as a response to stimuli or synaptic plasticity, naturally invites the use of  persistent homology associated to various filtrations of the complexes. Thus the need arose for a flexible and efficient computational package that is capable of generating complexes out of very large directed graphs and computing homology and persistent homology within reasonable compute time and use of memory. 

The software package \textsc{Flagser} presented in this article was developed as part of an EPSRC funded project, \emph{Topological Analysis of Neural Systems}. It is designed specifically for dealing with the kind of constructions and high performance computations required in this project, but is potentially applicable in a variety of other areas. The homology computation part of \textsc{Flagser} is based on a rather new addition to the topological data analysis software collection,  \textsc{Ripser} by Ulrich Bauer \cite{Bauer}. However, unlike \textsc{Ripser}, the \textsc{Flagser} package includes the capability to generate the (directed or undirected) flag complex of a finite (directed or undirected) graph from its connectivity matrix. In addition the reduction algorithm and coboundary computation methods were modified and optimized for use with  very large datasets. We also added a  flexible way of associating weights to vertices and/or edges in a given graph, and the capability to have user defined filtrations on the complex, resulting from the given weights. Another very useful added feature in \textsc{Flagser} is an Approximate function that allows the homology computation part of the program to compute Betti numbers and persistent homology within a user defined accuracy, and as a result dramatically reduce compute time. \textsc{Flagser} also contains the capability of computing homology with coefficients in odd characteristics. With \textsc{Flagser} we are able to perform the construction of the directed flag complex of the neocortical column of a rat (Link \ref{BBP-Rat}) on a laptop within less than 30 seconds, and compute a good approximation for it's homology within a few of hours on an HPC cluster (Figure \ref{resultstable}). Such computations were completely out of reach before. 

The usefulness of \textsc{Flagser} is clearly not limited to neuroscience. Topological analysis of graphs and complexes arising in both theoretical problems and applicable studies have become quite prevalent in recent years \cite{Farber, Brodzki}. Hence \textsc{Flagser}, in its original or modified form, could become a major computational aid in many fields. 

In the following, we discuss the main design choices and their impact on computations.
The last section contains real-world examples demonstrating the usefulness of the software package.  \textsc{Flagser} is an open source package, licensed under the LGPL 3.0 scheme. Some task specific modifications are available as well (See Section \ref{Source}).

\section{Mathematical background}
\label{background}

In this section, we give a brief introduction to  directed flag complexes of directed graphs, and to persistent homology.
For more detail, see \cite{Munkres}, \cite{Hatcher} and \cite{Reimann17}.

\begin{Defi}\label{Def-Graph}
  A \hadgesh{graph} $G$ is a pair $(V, E)$, where $V$ is a finite set referred to as the   \hadgesh{vertices} of $G$, and  $E$ is a subset of the set of unordered pairs $e=\{v, w\}$ of distinct points in $V$, which we call the \hadgesh{edges} of $G$. Geometrically the pair $\{u, v\}$ indicates that the vertices $u$ and $v$ are adjacent to each other in $G$.
  A \hadgesh{directed graph}, or a \hadgesh{digraph}, is  similarly a pair $(V,E)$ of vertices $V$ and edges $E$, except the edges are ordered pairs of distinct vertices, i.e.\ the pair $(v, w)$ indicates that there is an edge \hadgesh{from $v$ to $w$} in $G$.
 In a digraph we allow reciprocal edges, i.e. both $(u,v)$ and $(v,u)$ may be edges in $G$, but we exclude loops, i.e. edges of the form $(v,v)$.
\end{Defi}

Next we define the main topological objects considered in this article.

\begin{Defi}\label{Def-Ordered-SCx}
An \hadgesh{abstract simplicial complex on a  vertex set $V$} is a collection $X$ of non-empty finite  subsets  $\sigma\subseteq V$ that is closed under taking non-empty subsets. An \hadgesh{ordered simplicial complex on a vertex set $V$} is a collection of non-empty finite \hadgesh{ordered} subsets $\sigma\subseteq V$ that is closed under taking non-empty ordered subsets. Notice that the vertex set $V$ does not have to have any underlying order.
  In both cases the subsets $\sigma$ are called the \hadgesh{simplices of $X$}, and $\sigma$ is said to be a \hadgesh{$k$-simplex}, or a \hadgesh{$k$-dimensional} simplex, if it consists of $k+1$ vertices. If $\sigma \in X$ is a simplex and $\tau\subseteq \sigma$, then $\tau$ is said to be a \hadgesh{face} of $\sigma$. If $\tau$ is a subset of $V$ (ordered if $X$ is an ordered simplicial complex) that contains $\sigma$, then $\tau$ is called a \hadgesh{coface of $\sigma$}.
  If $\sigma = (v_0, \ldots, v_k)$ is a $k$-simplex in $X$, then the $i$-th face of $\sigma$ is the subset $\partial_i(\sigma) = (v_0, \ldots, \hat{v}_i,\ldots, v_k)$ where by $\hat{v}_i$ we mean omit the $i$-th vertex. For $d\ge 0$ denote by $X_d$ the set of all $d$-dimensional simplices of $X$.
\end{Defi}

Notice that if $X$ is any abstract simplicial complex with a vertex set $V$, then by fixing a linear order on the set $V$, we may think of each one of the simplices of $X$ as an ordered set, where the order is induced from the ordering on $V$. Hence any abstract simplicial complex together with an ordering on its vertex set determines a unique ordered simplicial complex. Clearly any two orderings on the vertex set of the same complex yield two distinct but isomorphic ordered simplicial complexes. Hence from now on we may assume that any simplicial complex we consider is ordered.

\begin{Defi}\label{Def-Directed-Flag-Complex}
  Let $G = (V,E)$ be a \hadgesh{directed} graph. The \hadgesh{directed flag complex} $\DFlag(G)$ is defined to be the ordered simplicial complex whose $k$-simplices are all  \hadgesh{ordered $(k+1)$-cliques}, i.e. $(k+1)$-tuples $\sigma = (v_0, v_1, \ldots, v_k)$, such that $v_i\in V$ for all $i$, and $(v_i,v_j)\in E$ for $i<j$.
  The vertex $v_0$ is called the \hadgesh{initial vertex} of $\sigma$ or the \hadgesh{source} of $\sigma$, while the vertex $v_k$ is called the \hadgesh{terminal vertex} of $\sigma$ or  the \hadgesh{sink} of $\sigma$.
\end{Defi}

Let $\F$ be a field, and let $X$ be an ordered simplicial complex. For any $k\ge 0$ let  $C_k(X)$  be the $\F$-vector space with basis given  by the set $X_k$ of all $k$-simplices of $X$. If $\sigma = (v_0, \ldots, v_k)$ is a $k$-simplex in $X$, then the $i$-th face of $\sigma$ is the subset $\partial_i(\sigma) = (v_0, \ldots, \hat{v}_i,\ldots, v_k)$, where $\hat{v}_i$ means that the vertex $v_i$ is removed to obtain the corresponding $(k-1)$-simplex. For $k\geq 1$, define the \hadgesh{$k$-th boundary map} $d_k\colon C_k(X) \to C_{k-1}(X)$ by
\[
  d_k(\sigma) = \sum_{i=0}^{k} (-1)^i \partial_i(\sigma) \in C_{k-1}(X).
\]
The boundary map has the property that  $d_{k}\circ d_{k+1}=0$ for each $k\geq 1$ \cite{Hatcher}. Thus $\im(d_{k+1})\subseteq \kernel(d_k)$ and one can make the following fundamental definition:

\begin{Defi}
  Let $X$ be an ordered simplicial complex. The $k$-th homology group of $X$ (over $\F$) denoted by $H_k(X)=H_k(X;\F)$ is defined to be
  \[
    H_k(X) = \kernel(d_k) / \im(d_{k+1}).
  \]
\end{Defi}

If $X$ is an ordered simplicial complex, and $Y\subseteq X$ is a subcomplex, then the inclusion $\iota\colon Y\to X$ induces a \hadgesh{chain map} $\iota_*\colon C_*(Y)\to C_*(X)$, i.e. a homomorphism $\iota_k\colon C_k(Y)\to C_k(X)$ for each $k\geq 0$, such that $d_k\circ\iota_k = \iota_{k-1}\circ d_k$. Hence one obtains an induced homomorphism for each $k\geq 0$,
\[\iota_*\colon H_k(Y) \to H_k(X).\]
This is a particular case of a much more general property of homology, namely its functoriality with respect to so called simplicial maps. We will not however require this level of generality in the discussion below. 

\begin{Defi}\label{Def-Filtered-SCx}
  A \hadgesh{filtered simplicial complex} is a simplicial complex $X$ together with a \hadgesh{filter function} $f\colon X\to \R$, i.e. a function satisfying the property that if $\sigma$ is a face of a simplex $\tau$ in $X$, then $f(\sigma) \le f(\tau)$.
  Given  $b\in\R\sqcup \{\infty\}$,  define the \hadgesh{sublevel complex} $X[b]$ to be the inverse image $f^{-1}\left( (-\infty, b] \right)$.
\end{Defi}

Notice that the condition imposed on the function $f$  in Definition \ref{Def-Filtered-SCx} ensures that for each $b\in\R\sqcup\{\infty\}$, the sublevel complex $X[b]$ is a subcomplex of $X$, and if $b\le b'$ then $X[b]\subseteq X[b']$. By default we identify $X$ with $X[\infty]$. Thus a filter function $f$ as above defines an increasing family of subcomplexes $X[b]$ of $X$ that is  parametrized by the real numbers.

If the values of the filter function $f$ are a finite subset of the real numbers, as will always be the case when working with actual data, then by enumerating those numbers by their natural order  we obtain a finite sequence of subspace inclusions
\[ X[1] \subseteq X[2] \subseteq \cdots \subseteq X[n] = X. \]
For such a sequence one defines its \hadgesh{persistent homology} as follows: for each $i\le j$ and  each $k\geq 0$ one has a homomorphism
\[h_k^{i,j}\colon H_k(X[i]) \to H_k(X[j]),\]
induced by the inclusion $X[i]\to X[j]$. The image of this homomorphism is a $k$-th persistent homology group of $X$ with respect to the given filtration, as it represents all homology classes that are present in $H_k(X[i])$ and carry over to $H_k(X[j])$. To use the common language in the subject, $\im(h_k^{i,j})$ consists of the classes that were ``born'' at or before filtration level $i$ and are still ``alive'' at filtration level $j$. The rank of $\im(h_k^{i,j})$ is a \hadgesh{persistent Betti number} for $X$ and is denoted by $\beta_k^{i,j}$.

Persistent homology can be represented as a collection of \hadgesh{persistence pairs}, i.e., pairs $(i,j)$, where $0\le i<j\le n$. The \hadgesh{multiplicity} of a persistence pair $(i,j)$ in dimension $k$ is, roughly speaking, the number of linearly independent $k$-dimensional homology classes that are born in filtration $i$ and die in filtration $j$. More formally, 
\[\mu_k^{i,j} = (\beta^{i,j-1}_k - \beta^{i-1,j-1}_k) - (\beta^{i,j}_k - \beta^{i-1,j}_k).\]
Persistence pairs can be encoded in \hadgesh{persistence diagrams} or \hadgesh{persistence barcodes}. See \cite{Edelsbrunner-Harer} for more detail.

\textsc{Flagser} computes the persistent homology of filtered directed flag complexes. The filtration of the flag complexes is based on filtrations of the 0-simplices and/or the 1-simplices, from which the filtration values of the higher dimensional simplices can then be computed by various formulas. When using the trivial filtration algorithm of assigning the value 0 to each simplex, for example, \textsc{Flagser} simply computes the ordinary simplicial Betti numbers. The software uses a variety of tricks to shorten the computation time (see below for more details), and outperforms many other software packages, even for the case of unfiltered complexes (i.e.\ complexes with the trivial filtration). In the following, we mainly focus on (persistent) homology with coefficients in finite fields, and if we don't mention it explicitly we consider the field with two elements $\mathbb{F}_2$.

%%%%%%%%%%%%%%%%%%%%%%%%%%%%%%%%%%%%%%%%%
%%%%%%%%%%%%%%%%%%%%%%%%%%%%%%%%%%%%%%%%%
%%%%%%%%%%%%%%%%%%%%%%%%%%%%%%%%%%%%%%%%%
%%%%%%%%%% SECTION 1.          %%%%%%%%%%%%%%%%%%%%%
%%%%%%%%%%%%%%%%%%%%%%%%%%%%%%%%%%%%%%%%%
%%%%%%%%%%%%%%%%%%%%%%%%%%%%%%%%%%%%%%%%%
%%%%%%%%%%%%%%%%%%%%%%%%%%%%%%%%%%%%%%%%%

\section{Representing the directed flag complex in memory}
An important feature of \textsc{Flagser} is its memory efficiency which allows it to carry out rather large computations on a laptop with a mere 16GB of RAM (see controlled examples in Section 6). In order to achieve this, the directed flag complex of a graph is generated on the fly and not stored in memory.
An $n$-simplex of the directed flag complex $\DFlag(G)$ corresponds to a directed $(n+1)$-clique in the graph $G = (V,E)$. Hence each $n$-simplex is identified with its list of vertices ordered by their in-degrees in the corresponding directed clique from $0$ to $n$. The construction of the directed flag complex is performed by  a double iteration:
\begin{enumerate}[(I)]
\item For each vertex of the graph we iterate over all simplices which have this vertex as the initial vertex.\label{i1}
\item Iterate \eqref{i1}  over all vertices. \label{i2}
\end{enumerate}
This procedure enumerates all simplices exactly once.

In more detail, \modifprev{consider a} fixed vertex $v_0$\modifprev{. The first iteration of Step \eqref{i1} starts}  by finding all vertices $w$ such that there exists a directed edge from $v_0$ to $w$. This \modifprev{produces} the set $S_1^{v_0}$ of all 1-simplices of the form $(v_0, w)$, i.e. all 1-simplices that have $v_0$ as their initial vertex.
\modifprev{In the second iteration of Step \eqref{i1}, construct for each 1-simplex $(v_0,v_1)\in S_1^{v_0}$ obtained in the first iteration, the set $S_2^{v_0,v_1}$  by finding} the intersection of the two sets of vertices that can be reached from $v_0$ and $v_1$, respectively. In other words, 
\[S_2^{v_0,v_1} = \{w \in V \;|\; (v_0, w)\in\ S_1^{v_0},\;\text{and}\; (v_1,w)\in S_1^{v_1}\}.\]
 This is the set of  all  2-simplices in $\DFlag(G)$ of the form $(v_0,v_1,w)$. 
\modifprev{Proceed by iterating Step \eqref{i1}, thus creating the prefix tree of the set of simplices with initial vertex $v_0$. The process must terminate, since the graph $G$ is assumed to be finite.} In Step \ref{i2} the same procedure is iterated over all vertices in the graph \modifprev{(but see \eqref{Ac} below.)}
\medskip

\modifprev{This algorithm has three main advantages:}
\begin{enumerate}[(a)]
  \item The construction of  the \modifprev{directed flag complex is very highly parallelizable. In theory one could use one CPU for each of the vertices of the graph, so that each CPU computes only  the simplices with this vertex as initial vertex.} \label{Aa}
  
  \item Only the graph is loaded into memory. The directed flag complex, which is usually much bigger, does not have any impact on the memory footprint. \label{Ab}
  
  \item  \modifprev{The procedure skips branches of the iteration based on the prefix. The idea is that if a simplex is not contained in the complex in question, then no simplex containing it as a face can be in the complex. Therefore, if we computed that a simplex $(v_0, \ldots, v_k)$ is not contained in our complex, then we don't  have to iterate Step \eqref{i1} on its vertices to compute the vertices $w$ that are reachable from each of the $v_i$. This allows us to  skip the full iteration branch of all simplices with prefix $(v_0, \ldots, v_k)$. This is particularly useful for iterating over the simplices in a subcomplex.}  \label{Ac}
\end{enumerate}

In the code, the adjacency matrix (and for improved performance its transpose) are both stored as bitsets\footnote{Bitsets store a single bit of information (0 or 1) as one physical bit in memory, whereas usually every bit is stored as one byte (8 bit). This representation improves the memory usage by a factor of 8, and additionally allows the compiler to use much more efficient bitwise logical operations, like ``AND'' or ``OR''. Example of the bitwise logical ``AND'' operation: \texttt{0111 \& 1101 = 0101}.}. \modifprev{This makes the computation of the directed flag complex more efficient. Given a $k$-simplex $(v_0, \ldots, v_k)$ in $\DFlag(G)$, the set of vertices
\[\{w\in V\;|\; (v_0, \ldots, v_k, w)\in \DFlag(G)\}\]
is computed as} the intersection of the sets of vertices that the vertices $v_0, \ldots, v_k$ have a connection to.
Each such set is given by the positions of the \texttt{1}'s in the bitsets of the adjacency matrix in the rows corresponding to the vertices $v_0, \ldots, v_k$, and thus the set intersection can be computed as the logical ``AND'' of those bitsets. 

\vspace{1em}

In some situations, we might have enough memory available to load the full complex into memory.
For these cases, there exists a modified version of the software that stores the tree described above in memory, allowing to
\begin{itemize}
  \item \modifprev{perform computations on the complex such as counting, or computing coboundaries etc.,} without having to recompute the intersections described above, and 
  
  \item associate data to each simplex with fast lookup time.
\end{itemize}
The second part is useful for assigning unique \modifprev{identifiers} to each simplex. Without the complex in memory a hash map with custom hashing function has to be created for this purpose, which slows down the lookup of the global \modifprev{identifier} of a simplex dramatically.

%%%%%%%%%%%%%%%%%%%%%%%%%%%%%%%%%%%%%%%%%
%%%%%%%%%%%%%%%%%%%%%%%%%%%%%%%%%%%%%%%%%
%%%%%%%%%%%%%%%%%%%%%%%%%%%%%%%%%%%%%%%%%
%%%%%%%%%% SECTION 4           %%%%%%%%%%%%%%%%%%%%%
%%%%%%%%%%%%%%%%%%%%%%%%%%%%%%%%%%%%%%%%%
%%%%%%%%%%%%%%%%%%%%%%%%%%%%%%%%%%%%%%%%%
%%%%%%%%%%%%%%%%%%%%%%%%%%%%%%%%%%%%%%%%%

\section{Computing the coboundaries}

In \textsc{Ripser}, \cite{Bauer}, cohomology is computed instead of homology. The reason for this choice is that for typical examples, the coboundary matrix contains many more columns that can be skipped in the reduction algorithm than the boundary matrix \cite{BauerPres}.

\modifprev{If $\sigma$ is a simplex in a complex $X$, then $\sigma$ is a canonical basis element in the chain complex for $X$ with coefficients in the field $\F$. The hom-dual $\sigma^*$, i.e. the function that associates $1\in\F$ with $\sigma$ and 0 with any other simplex, is a basis element in the cochain complex computing the cohomology of $X$ with coefficients in $\F$.  The differential in the cochain complex requires considering the coboundary simplices of any given simplex, i.e. all those simplices that contain the original simplex as a face of codimension 1.}

In \textsc{Ripser}, the coboundary simplices are never stored \modifprev{in memory.  Instead, the software enumerates the coboundary simplices of any given simplex by a clever indexing technique, omitting all coboundary simplices that are not contained in the complex defined by the  filtration stage being computed.
However, this technique only works if the number of theoretically possible coboundary entries is sufficiently small.}
In the \hadgesh{directed} situation, and when working with a directed graph with a large number of vertices, this number is typically too large to be enumerated by a 64-bit integer.
Therefore, in \textsc{Flagser} we precompute the coboundary simplices into a sparse matrix, achieving fast coboundary iteration at the expense of  longer initialization time and more memory load.
An advantage of this technique is that the computation of the \modifprev{coboundary matrices} can be fully parallelized, again splitting the simplices into groups indexed by their initial vertex.

Given an ordered list of vertex identifiers $(v_0, \ldots, v_k)$ characterising a directed simplex \modifprev{$\sigma$, the list of its coboundary simplices} is computed by enumerating all possible ways of inserting a new vertex into this list \modifprev{to create a simplex in the coboundary of $\sigma$.} Given a position $0\le p\le k+1$, the set of vertices that could be inserted in that position is given by the \modifprev{intersection of the} set of vertices that have a connection from each of the vertices $v_0, \ldots, v_{p-1}$ and \modifprev{the set of vertices that have a connection} to each of the vertices $v_p, \ldots, v_k$.
This intersection can again be efficiently computed by bitwise logical operations of the rows of the adjacency matrix and its transpose (the transpose is necessary to efficiently compute the incoming connections of a given vertex). 

With the procedure described above, we can compute the coboundary of every $k$-simplex. The simplices in this representation of the coboundary are given by their ordered list of vertices, and if we want to turn this into a coboundary matrix we have to assign consecutive numbers to the simplices of the different dimensions. We do this by creating a hash map, giving all simplices of a fixed dimension consecutive numbers starting at 0. This numbering corresponds to choosing the basis of our chain complex by fixing the ordering, and using the hashmap we can lookup the number of each simplex in that ordering. This allows us to turn the coboundary information computed above into a matrix representation. 
If the flag complex is stored in memory, creating a hash map becomes unnecessary since we can use the complex as a hash map directly, adding the dimension-wise consecutive identifiers as additional information to each simplex.
This eliminates the lengthy process of hashing all simplices and is therefore useful in situations where memory is not an issue.

%%%%%%%%%%%%%%%%%%%%%%%%%%%%%%%%%%%%%%%%%
%%%%%%%%%%%%%%%%%%%%%%%%%%%%%%%%%%%%%%%%%
%%%%%%%%%%%%%%%%%%%%%%%%%%%%%%%%%%%%%%%%%
%%%%%%%%%% SECTION 5           %%%%%%%%%%%%%%%%%%%%%
%%%%%%%%%%%%%%%%%%%%%%%%%%%%%%%%%%%%%%%%%
%%%%%%%%%%%%%%%%%%%%%%%%%%%%%%%%%%%%%%%%%
%%%%%%%%%%%%%%%%%%%%%%%%%%%%%%%%%%%%%%%%%

\section{Performance considerations}
\modifprev{\textsc{Flagser} is optimized for large computations. In this section we describe some of the algorithmic methods by which this is carried out.}

\subsection{Sorting the columns of the coboundary matrix}
The reduction algorithm iteratively reduces the columns in order of their filtration value (i.e.\ the filtration value of the simplex that this column represents).
Each column is reduced by adding previous (already reduced) columns until the pivot element does not appear as the pivot of any previous column.
While reducing all columns the birth and death times for the persistence diagram can be read off: roughly speaking, if a column with a certain filtration value $f_0$ ends up with a pivot that is the same as the pivot of a column with a \emph{higher} filtration value $f_1$, then this column represents a class that is born at time $f_0$ and dies at time $f_1$. 
Experiments showed that the order of the columns is crucial for the performance of this reduction algorithm, so \textsc{Flagser} sorts the columns that have the same filtration values by an experimentally determined heuristic (see below).
This is of course only useful when there are many columns with the same filtration values, for example when taking the trivial filtration where every simplex has the same value. In these cases, however, it can make very costly computations tractable.
\medskip

The heuristic to sort the columns is applied  as follows:
\begin{enumerate}
  \item Sort columns in ascending order by their filtration value. This is necessary for the algorithm to produce the correct results.
  
  \item  Sort columns in ascending order  by the number of non-zero entries in the column. This tries to keep the number of new non-zero entries that a reduction with such a column produces small.
  \item Sort columns in descending order by the position of the pivot element of the unreduced column. This means that the larger the pivot (i.e.\ the higher the index of the last non-trivial entry), the lower the index of the column is in the sorted list.  The idea behind this is that  ``long columns'', i.e. columns whose pivot element has a high index, should be as sparse as possible. Using such a column to reduce the columns to its right is thus more economical. Initial coboundary matrices are typically very sparse. Thus ordering the columns in this way gives more sparse columns with large pivots to be used in the reduction process.
  
  \item Sort columns in descending order by their ``pivot gap'', i.e. by the distance between the pivot and the next nontrivial entry with a smaller row index.  A large pivot gap is desirable in case the column is used in reduction of columns to its right,   because when using such a column to reduce another column, the added  non-trivial entry in the column being reduced appears with a lower index. This may yield   fewer reduction steps on average. 
\end{enumerate}

\begin{Rma}
  We tried various heuristics---amongst them different orderings of the above sortings---but the setting above produced by far the best results. Running a genetic algorithm in order to find good combinations and orderings of sorting criteria produced slightly better results on the small examples we trained on, but it did not generalize to bigger examples. It is known that finding the perfect ordering for this type of reduction algorithm is NP-hard \cite{Yannakakis}.
\end{Rma}

\subsection{Approximate computation}
Another \modifprev{method} that speeds up computations considerably is \modifprev{the use of} approximate computation. \modifprev{This simply means allowing \textsc{Flagser}} to skip columns that need  many reduction steps.
The reduction of these columns takes the most time, so skipping even only the columns with extremely long reduction chains gives significant performance gains.
By skipping a column, the rank of the matrix can be changed by at most one, so the error can be explicitly bounded.
This is especially useful for computations with trivial (i.e.\ constant) filtrations, as the error in the resulting Betti numbers is easier to interpret than the error in the set of persistence pairs.

Experiments on smaller examples show that the theoretical error that can be computed from the number of skipped columns is usually much bigger than the actual error. Therefore, it is usually possible to quickly get rather reliable results that afterwards can be refined with more computation time.
See \autoref{fig:approximate_computation} for an example of the speed-up gained by approximate computation. The user of \textsc{Flagser} can specify an approximation level, which is given by a number (defaulting to infinity). This number determines after how many reduction steps the current column is skipped. The lower the number, the more the performance gain (and the more uncertainty about the result). At the current time it is not possible to specify the actual error margin for the computation, as it is very complicated to predict the number of reductions that will be needed for the different columns.

\begin{figure}[htpb]
  \centering
  \begin{tabular}{|| l || c | c| c | c | c ||}
    \hline
    Approximation & none & $10^4$ & $10^3$ & $10^2$ & $10^1$ \\
    \hline
    Theoretical error & 0 & 322 & 1380 & 7849 & 192638   \\
    Actual error & 0 & 0 & 19 & 58 & 65   \\
    Computation time & 113s & 49s & 43s & 41s & 39s   \\
    \hline
\end{tabular}
  \caption{Theoretical and actual error in the computation of the first Betti number by skipping the reduction of columns that require certain numbers of reduction steps. The cell counts are 4656, 35587 and 1485099 in dimensions $0, 1$ and $2$ \modifprev{respectively}, and the first Betti number is 3321. The computations were made as averages over ten runs on a MacBook Pro, 2.5 GHz Intel Core i7.}
  \label{fig:approximate_computation}
\end{figure}

\subsection{Dynamic priority queue}
When reducing the matrix, \textsc{Ripser} stores the column that is currently reduced as a priority queue, ordering the entries in descending order by filtration value and in ascending order by an arbitrarily assigned index for entries with the same filtration value. Each index can appear multiple times in this queue, but due to the sorting all repetitions are next to each other in a contiguous block. To determine the final entries of the column after the reduction process, one has to loop over the whole queue and sum the coefficients of all repeated entries.  

This design makes the reduction of columns with few reduction steps very fast but gets very slow (and memory-intensive) for columns with a lot of reduction steps:
each time we add an entry to the queue it has to ``bubble'' through the queue and find its right place.
Therefore, \textsc{Flagser} enhances this queue by dynamically switching to a hash-based approach if the queue gets too big:
after a certain number of elements were inserted into the queue, the queue object starts to track the ids that were inserted in a hash map, storing each new entry with coefficient 1.
If a new entry is then already found in the hash map, it is not again inserted into the queue but rather the coefficient of the hash map is updated, preventing the sorting issue described above.
When removing elements from the front, the queue object then takes for each element that is removed from the front of the queue the coefficient stored in the hash map into account when computing the final coefficient of that index. 

\subsection{Apparent pairs}
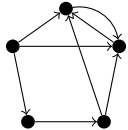
\begin{figure}\begin{tikzpicture}[scale=1]\def\s{0.5}
\node[circle, fill=black,scale=\s] (0) at (0,0){};
\node[circle, fill=black,scale=\s] (1) at (-.2,1){};
\node[circle, fill=black,scale=\s] (2) at (.5,1.5){};
\node[circle, fill=black,scale=\s] (3) at (1.2,1){};
\node[circle, fill=black,scale=\s] (4) at (1,0){};
\draw[->] (0) -- (4);
\draw[->] (1) -- (0);
\draw[->] (1) -- (2);
\draw[->] (1) -- (3);
\draw[->] (2) to [bend left=45] (3);
\draw[->] (3) -- (2);
\draw[->] (4) -- (2);
\draw[->] (4) -- (3);
\end{tikzpicture}
\caption{A digraph for which using apparent pairs to compute the homology of the directed flag complex gives the wrong results.}\label{fig:APgraph}\end{figure}

In \textsc{Ripser}, so-called \hadgesh{apparent pairs} are used in order to skip some computations completely.
These rules are based on discrete Morse theory \cite{Forman98}, but they only apply for \hadgesh{non-directed} simplicial complexes.
Since we consider directed graphs, the associated flag complex is  a semi-simplicial rather than simplicial complex, so we cannot use this simplification.
Indeed, experiments showed that enabling the skipping of apparent pairs yields wrong results for directed flag complexes of certain directed graphs.
For example, when applied to the directed flag complex of the graph in Figure~\ref{fig:APgraph}, using apparent pairs gives $\beta_1=\beta_2=0$ and not using apparent pairs gives $\beta_1=\beta_2=1$, which is visibly the correct answer.

Experiments showed that when applied to directed flag complexes that arise from graphs without reciprocal connections (so that the resulting  directed flag complex is a simplicial complex), enabling the skipping of apparent pairs does not give a big performance advantage. In fact, we found a significant reduction in run time for subgraphs of the connectome of the reconstruction of the neocortical column based data averaged across the five rats that were used in \cite{Reimann17}, the Google plus network and the Twitter network from KONECT (Link \ref{KONECT}).

%%%%%%%%%%%%%%%%%%%%%%%%%%%%%%%%%%%%%%%%%
%%%%%%%%%%%%%%%%%%%%%%%%%%%%%%%%%%%%%%%%%
%%%%%%%%%%%%%%%%%%%%%%%%%%%%%%%%%%%%%%%%%
%%%%%%%%%% SECTION 4.          %%%%%%%%%%%%%%%%%%%%%
%%%%%%%%%%%%%%%%%%%%%%%%%%%%%%%%%%%%%%%%%
%%%%%%%%%%%%%%%%%%%%%%%%%%%%%%%%%%%%%%%%%
%%%%%%%%%%%%%%%%%%%%%%%%%%%%%%%%%%%%%%%%%

\section{Sample Computations}

In this section we present some sample computations carried out using \textsc{Flagser}. Almost all computations were done on a laptop PC with 32GB of RAM and a 2-core processor. The data files were either taken from online public domain sources or generated for the purpose of computation, as detailed below.

\begin{figure}[h!]
{\small
\begin{tabular}{| l || l | l | l | l | l | l | l |}
\hline
Dataset & Vertices & Edges & Density \%&\texttt{Count} & \texttt{Homology} \\
\hline\hline
BE & 43 & 336 &18.172& \phantom{0}0.01s, \phantom{56}3.32MB & \phantom{10h44m}0.08s, \phantom{56}3.81MB\\
GP & 23.6K & 39.2K &\phantom{0}0.007& \phantom{0}0.50s, 564.40MB & \phantom{10h44m}2.09s, 569.54MB\\
IN & 410 & 2765 &\phantom{0}1.644& \phantom{0}0.04s, \phantom{56}3.68MB & \phantom{10h44m}1.27s, \phantom{56}9.13MB\\
JA$^\ast$ & 198 & 2742 &\phantom{0}6.994& 49.63s, \phantom{56}3.57MB & 10h44m3.00s, \phantom{0}75.79GB\\
MQ & 62 & 1187 &30.879& \phantom{0}0.09s, \phantom{56}3.44MB & \phantom{10h44m}9.41s, \phantom{5}50.22MB \\
PR & 2239 & 6452 &\phantom{0}0.129& \phantom{0}0.04s, \phantom{56}8.63MB & \phantom{10h44m}0.07s, \phantom{56}9.17MB\\
\hline
FN & 40 & 760 &47.500& 2m23.26s, 3.46MB & /\\
BB$^{\dagger\ast}$ & 31.3K & 7.8M &\phantom{0}0.793& \phantom{2m}23.76s, 1.08GB & 12h15m3.00s, 52.96 GB\\
BA & 280 & 1911 &\phantom{0}2.437& \phantom{2m0}0.00s, 3.52MB & \phantom{12h15m}0.03s, \phantom{0}3.84MB\\
CE & 279 & 2194 &\phantom{0}2.819& \phantom{2m0}0.02s, 3.55MB & \phantom{12h15m}0.12s, \phantom{0}4.55KB\\
ER$^{\dagger\ast}$ & 31.3K& 7.9M &\phantom{0}0.804& \phantom{2m}15.06s, 1.08GB & 1h30m51.00s, \phantom{0}4.98GB\\
\hline
\end{tabular}
}
\caption{Performance of \textsc{Flagser} on a variety of datasets: BE = Windsurfers, GP = Google+, IN = Infectious, JA = Jazz, MQ = Macaques, PR = Protein, FN = 40cycle19, BB = Blue Brain Project neocortical column, BA = Barabasi-Albert Graph, CE = C-Elegans, ER = Erd\H{o}s-R\'enyi Graph.}
\label{resultstable}
\end{figure}

In the table above, the dagger symbol on the dataset abbreviation means that the approximate function of \textsc{Flagser} was used with approximation parameter $100000$ (the code stops reducing a column after 100000 steps). This allows for computation of homology in a reasonable compute time. The asterisk stands for a homology computation that was made on a single core of an HPC cluster. This was done because the memory requirements of the homology computation were too large to perform it on a laptop.

The first six datasets were taken from the public domain collection KONECT (the Koblenz Network Collection -- Link \ref{KONECT}). We chose the following datasets: HumanContact/Windsurfers, Social/Google+, HumanContact/Infectious, HumanSocial/Jazz musicians, Animal/Macaques, Metabolic/Human protein (Figeys). See Figure \ref{resultstable} for detail on performance. Two of the networks came as  weighted graphs, which allowed us to compute persistent homology in those two cases (Figure \ref{persistence}). The filtration was determined by assigning to each simplex  the maximal weight of an edge it contains.

\begin{figure}[ht!]
{\small
\begin{tabular}{| l || l |}
\hline
Dataset &\texttt{Persistence}\\
\hline\hline
BE\ & \phantom{0}0.09s, \hskip.2in \phantom{0}3.91MB \\
MQ & 11.47s, \hskip.2in 69.54MB \\
\hline
\end{tabular}
}
\caption{Performance of \textsc{Flagser} computing persistence.}
\label{persistence}
\end{figure}

The last five datasets are: 19th step of the Rips filtration of the 40-cycle graph (equipped with the graph metric), a sample of a Blue Brain Project reconstruction of the neocortical column or a rat \cite{Reimann17}, an Erd\H{o}s-R\'enyi graph with a similar number of vertices and connection probability as such a network, a C. Elegans brain network \cite{Varshney} (Link \ref{C-Elegans}) and a Barab\'asi-Albert graph with a similar number of vertices and connection probability as a C. Elegans brain network.

Next, we tested the performance of \textsc{Flagser}  with respect to several values of the Approximate parameter. The performed the computation on two graphs: a)  a Blue Brain Project neocortical column reconstruction, denoted above by BB (Figure \ref{performance}), and an Erd\H{o}s-R\'enyi graph with the same number of vertices and connection probability (Figure \ref{ERperformance}). These computations were carried out on a single core of an HPC cluster.   

\hbox{}

\begin{figure}[ht!]
{\small
\begin{tabular}{| l || l | l || l | l | l | l | l | l | }
\hline
Approx.& Time & Memory & $\beta_0$ &$\beta_1$  & $\beta_2$ & $\beta_3$ & $\beta_4$ & $\beta_5$\\
\hline\hline
1 & 40m44s & 11.77GB & 1 & 14,640 & 9,831,130 & 30,219,098 & 513,675 &  785\\
10 & 38m44s & 11.96GB & 1 &  15,148&  14,831,585 & 12,563,309 & 13,278 & 40\\
100 & 40m29s & 12.69GB & 1 &  15,188 &  13,956,297 & 7,171,621 & 9,898 & 37\\
1000 & 55m40s & 15.42GB & 1 & 15,283 & 14,598,249 & 5,780,569 & 9,822 & 37\\
10000 & 2h24m17s & 23.78GB & 1 & 15,410 &  14,872,057 & 5,219,666 & 9,821 & 37\\
100000 & 12h15m3s & 52.96GB & 1&  15,438&  14,992,658 & 4,951,674 & 9,821 & 37\\
\hline
\end{tabular}
}
\caption{Performance of \textsc{Flagser} on the Blue Brain Project graph with respect to different values of the Approximate parameter. The Betti numbers computed are an approximation to the true Betti numbers. See Figure \ref{skipped} for accuracy.}
\label{performance}
\end{figure}

\begin{figure}[ht!]
{\small
\begin{tabular}{| l ||  r  | r |  r |  r  | r  |}
\hline
Approx. & $\delta_1$ & $\delta_2$ & $\delta_3$ & $\delta_4$ & $\delta_5$ \\
\hline
1 & 		7,724,934	&56,746,602 	&19,092,109 	&664,481 	&4,162\\
10 & 		6,714,811	&15,944,792 	&796,559 		&372 	&0\\
100 & 	179,794 	&4,107,019	&6,414 		&0 		&0\\
1000 &	14,209  	&1,902,229 	&42 			&0 		&0\\
10000 & 	1,002 	&1,054,397 	&0 			&0 		&0\\
100000 &  27		&664,857 		&0 			&0 		&0\\
\hline
\end{tabular}
}
\caption{Number of skipped columns of the coboundary matrix for the Blue Brain Project graph with respect to different values of the Approximate parameter. The approximation accuracy of $\beta_i$ depends theoretically on the number of columns skipped in the coboundary matrices for $\delta_{i-1}$ and $\delta_i$.}
\label{skipped}
\end{figure}

\begin{figure}[ht!]
\includegraphics[width=0.6\textwidth]{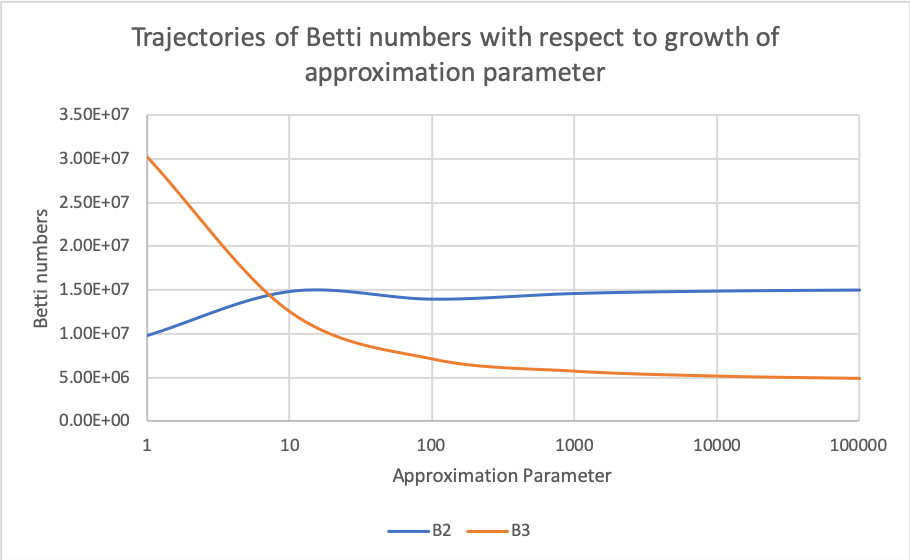}
\caption{Trajectories of $\beta_2$ and $\beta_3$ as a function of the approximation parameter. A change in order of magnitude in the approximation parameter results in a small change in the Betti number computed.}
\label{Convergence}
\end{figure}

Consider for example $\beta_4$. By Figures \ref{performance} and \ref{skipped} the true value is $\beta_4 = 9821$, since in that case no columns are skipped in $\delta_3$ and $\delta_4$. Thus perfect accuracy is achieved with Approximate parameter 10000. However, with Approximate parameters 1000 and 100 one has a theoretical error of 42 and 6414 respectively, whereas the actual error is 1 and 77 respectively. In the case of $\beta_2$ we do not know the actual answer. However, since the total number skipped with Approximate parameter 100000 is 664884, and the approximate value computed for $\beta_2$ is 14992658, the calculation is accurate within less than $\pm 5\%$. Similarly, $\beta_1$ is computed to be 15438 with 27 skipped columns, which stands for an accuracy at least as good as $\pm0.2\%$. For $\beta_3$ the theoretical bound is clearly not as good, but here as in all other cases where precise computation was not carried out, the trajectory of the approximated Betti numbers as the approximation parameter grows suggests that the number computed is actually much closer to the real Betti numbers than the theoretical error suggests (Figure \ref{Convergence}).

We performed a similar computation with an  Erd\H{o}s-R\'enyi  random directed graph with the same number of vertices as the Blue Brain Project column and the same average probability of connection (Figures \ref{ERperformance} and \ref{ERskipped}). Notice that there is a very substantial difference in performance between the Blue Brain Project graph and a random graph with the same size and density parameters. The difference can be explained by the much higher complexity of the Blue Brain Project graph, as witnessed by the Betti numbers.

\begin{figure}[ht!]
{\small
\begin{tabular}{| l || l | l |  l |  l |  l |  l | }
\hline
Approx. & Time & Memory & $\beta_0$ & $\beta_1$ & $\beta_2$ & $\beta_3$ \\
\hline\hline
1 		& 4m45s 		& 3.25GB 		& 1 & 19,675	 & 14,675,052 	& 40\\
10 		& 4m31s 		&3.10GB 		& 1 & 19,981	& 9,204,135 	& 4\\
100 		& 4m44s 		& 3.18GB 		& 1 & 19,981	& 8,074,149 	& 4\\
1000 	& 5m48s 		& 3.33GB 		& 1 &  19,981	& 7,921,730 	& 4\\
10000 	& 14m18s 	& 3.72GB 		& 1 &  19,982	& 7,876,858 	& 4\\
100000 	& 1h30m51s 	&  4.98GB		& 1 &  19,984 	& 7,857,917 	& 4\\
\hline
\end{tabular}
}
\caption{Performance of \textsc{Flagser} on an Erd\H{o}s-R\'enyi graph with respect to different values of the Approximate parameter. The Betti numbers range from dimension 0 to dimension 3.}
\label{ERperformance}
\end{figure}

\begin{figure}[ht!]
{\small
\begin{tabular}{| l || l | l | l |}
\hline
Approx. & $\delta_1$ & $\delta_2$ & $\delta_3$ \\
\hline\hline
1 		& 7,692,010 	& 719,041 	&124\\
10 		& 1,501,909  	& 3 			& 0\\
100 		& 371,920  	& 0 			& 0\\
1000 	&  219,501  	& 0  			& 0\\
10000 	&  174,628  	& 0 			& 0\\
100000 	& 155,685  	& 0  			& 0\\
\hline
\end{tabular}
}
\caption{Number of skipped columns of the coboundary matrix for an Erd\H{o}s-R\'enyi graph with respect to different values of the Approximate parameter in dimensions from 1 to 3.}
\label{ERskipped}
\end{figure}

\begin{figure}
\includegraphics[scale=.7]{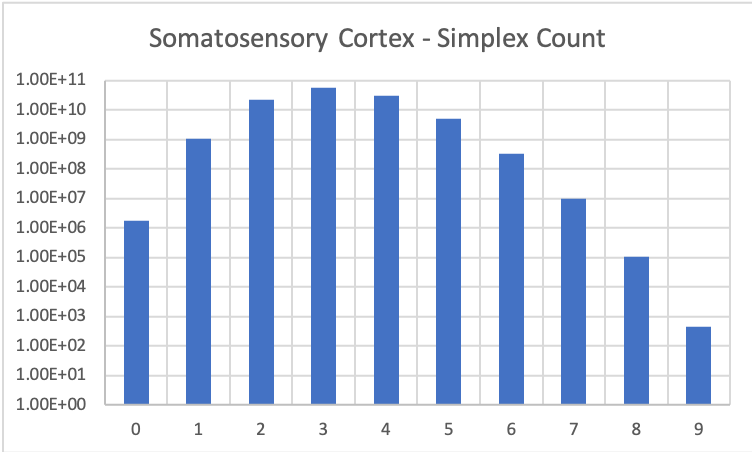}
\caption{Cell count in the Blue Brain reconstruction of the somatosensory cortex. Computation was run on an HPC cluster using 256 CPUs, required 55.69GB of memory and took 7.5 hours to complete. }
\label{Somatosensory}
\end{figure}

\begin{figure}
\includegraphics[scale=.5]{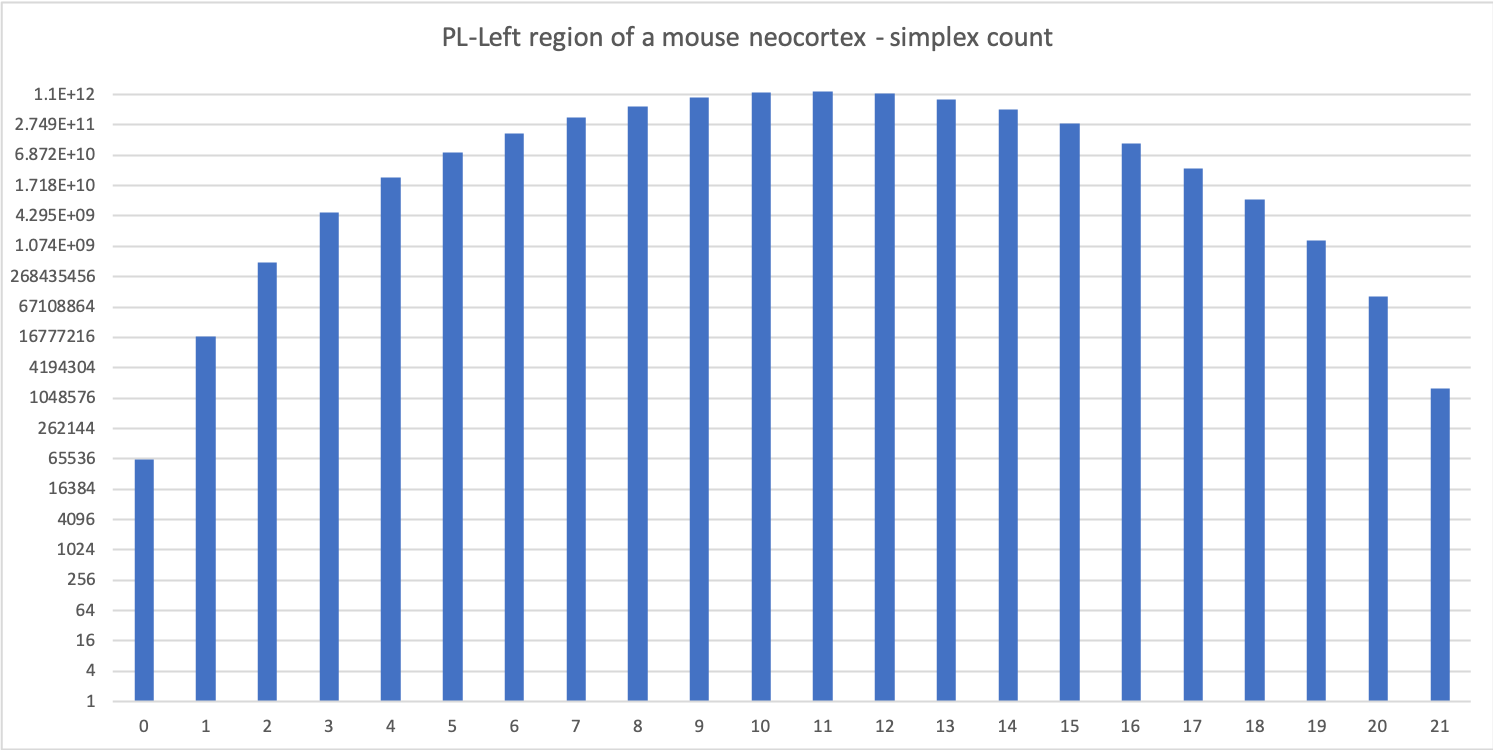}
\caption{The simplex count for the PL-Left region of a mouse neocortex (local connections only). The complex is 21-dimensional with more than 1.2 trillion 11-dimensional simplices. Computation was run in parallel on two nodes of an HPC cluster with 256 CPUs each, required 1GB of memory and took 5 days to complete (run as 10 different jobs of 24 hours each).}
\label{PL}
\end{figure}

\begin{figure}
\centering
\begin{subfigure}{0.4\textwidth}
\includegraphics[scale=.3]{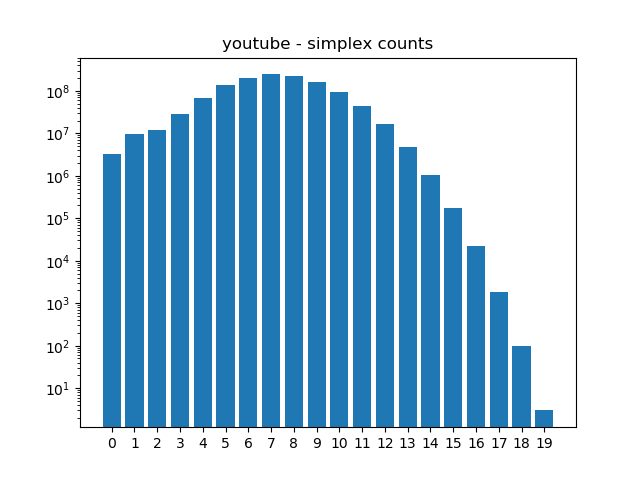}
\centering
\caption{Runtime: 48m45s, Memory: 2.5GB}
\end{subfigure}
\begin{subfigure}{0.4\textwidth}
\includegraphics[scale=.3]{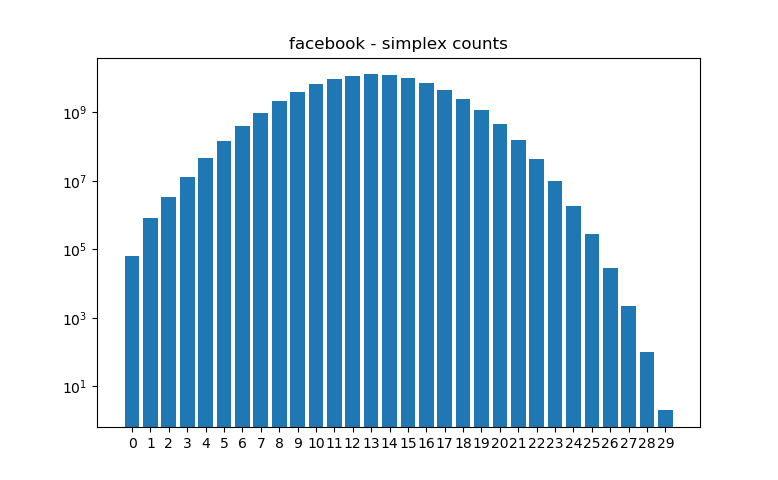}
\caption{Runtime: 6h24m50s, Memory: 40.5GB}
\end{subfigure}
\vspace{1em}
\begin{subfigure}{0.4\textwidth}
\centering
\includegraphics[scale=.3]{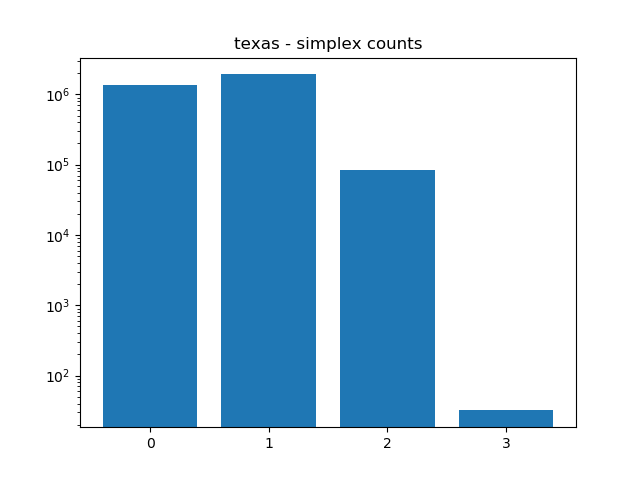}
\caption{Runtime: 3m49s, Memory: 1GB}
\end{subfigure}
\begin{subfigure}{0.4\textwidth}
\centering
\includegraphics[scale=.3]{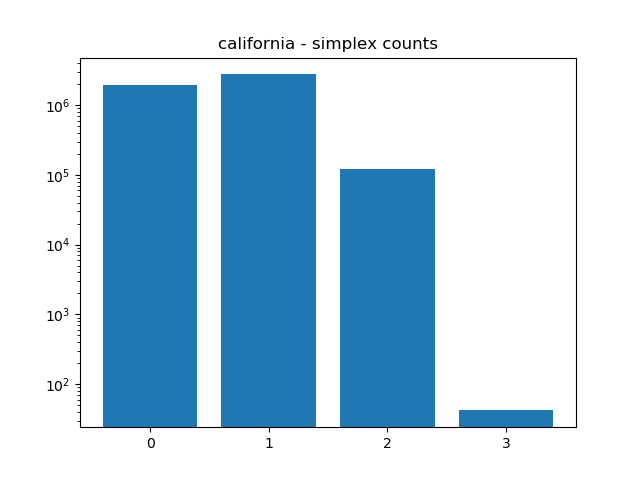}
\caption{Runtime: 5m19s, Memory: 1.5GB}
\end{subfigure}
\caption{Four large datasets taken from the Koblenz Network Collection (Link \ref{KONECT}). (a) Social network of YouTube users and their friendship connections, (b) Friendship data of Facebook users, (c) Road network of the State of Texas, a node is an intersection between roads or road endpoint and edges are road segments and (d) Road network of the State of California. All computations were carried out on an HPC cluster using 256 cores.}
\label{othercount}
\end{figure}

A remarkable feature of \textsc{Flagser} is its parallelisable construction algorithm of the directed flag complex. To demonstrate the capability of the software, we computed some large scale examples. In particular, we computed the cell counts of two sample circuits produced by the Blue Brain Project. The first is a reconstruction of the somatosensory cortex of a rat (Figure~\ref{Somatosensory}), and the second is the so called PL region of the reconstruction of the full neocortex of a mouse (Figure~\ref{PL}). We also computed the simplex counts on some large non-biological networks (Figure~\ref{othercount}).

The somatosensory cortex (Figure \ref{Somatosensory}) is a digraph with about 1.7 million vertices and 1.1 billion edges (average connection probability $0.04\%$).
The Blue Brain Project model of the full mouse cortex  is available online (Link \ref{BBP-Mouse}). The full reconstruction consists of approximately 10 million neurons. Its synaptic connections fall into three categories. Local (short distance) connections, intermediate distance connections and long distance connections. We computed the cell count of the PL region of the left hemisphere in the reconstruction, consisting of approximately 64,000 neurons and 17.5 million local connections (average connection probability $0.4\%$). We did not consider any connections but the local ones. 

Interestingly the number of neurons in the model of the PL region (Link \ref{BBP-Mouse}) is roughly twice that in the reconstruction of the neocortical column studied in \cite{Reimann17}, and the probability of connection is roughly half. Hence we expected a comparable performance of \textsc{Flagser}. The results are therefore quite surprising. The likely reason for the difference in computation speeds is the higher dimensionality of the model of neocortex of the mouse, resulting in a gigantic number of simplices in the middle dimensions (approximately 1.1 trillion 11-dimensional simplices), compared to the neocortical column of the rat, with the former containing simplices of dimension 21 (see Figure~\ref{PL}) and the latter only containing simplices up to dimension 7. 

These computation were done on an adapted version of \textsc{Flagser} available at \ref{Flag-Adapt]}). This alternate version has been adapted in the following ways: the code to compute homology is removed as this is not needed, the option to load the graph from compressed sparse matrix format is included as this speeds up the time required to read the graph, the graph can be stored using Google's dense hash map to reduce the memory needed for very large graphs, and partial results are printed during run time which allows for the computation to be resumed if it is stopped (for example, if the process is killed due to time limitations on HPC).

\section{Availability of source code}
\label{Source}

\textsc{Flagser}  was developed as a part of an EPSRC funded project on applications of topology in neuroscience. As such \textsc{Flagser} is an open source package, licensed under the LGPL 3.0 scheme. The code can be found on the GitHub page of Daniel L\"utgehetmann who owns the copyright for the code and will maintain the page (Link \ref{Flagser}).

In the course of using \textsc{Flagser} in the project three adapted versions were developed, and can be found in the GitHub page of Jason Smith (Link \ref{Flag-Adapt}), who will maintain the page.  The following modifications are available. 
\begin{itemize}
\item \textsc{Flagser-count}: A modification of a basic procedure in \textsc{Flagser}, optimised for very large networks.
\item \textsc{Deltser}:  A modification of  \textsc{Flagser} that computes persistent homology on any finite delta complex (like a simplicial complex, except more than one simplex may be supported on the same set of vertices).
\item  \textsc{Tournser}: The \hadgesh{flag tournaplex} of a directed graph is the delta complex whose simplices are all tournaments  that occur in the graph.  \textsc{Tournser} is an adaption of \textsc{Flagser} which computes (persistent) homology of the tournaplex of a finite directed graph. Tournser uses two filtrations that occur naturally in tournaplexes, but can also be used in conjunction with other filtrations. 
\end{itemize}

\section{Acknowledgements} The authors acknowledge support by an EPSRC grant EP/P025072/1 - \hadgesh{Topological Analysis of Neural Systems}. 

\section{Web Links}
\label{Links}

\subsection{\text[BBP]} \url{https://www.epfl.ch/research/domains/bluebrain/}{}.
\label{BBP}

\subsection{\text{[BBP-Mouse]}} 
\url{https://portal.bluebrain.epfl.ch/resources/models/mouse-projections/}{}. 
\label{BBP-Mouse}

\subsection{\text{[BBP-Rat]}} \url{https://bbp.epfl.ch/nmc-portal/downloads}{}.
\label{BBP-Rat}

\subsection{\text{[C-Elegans]}} \url{http://www.wormatlas.org/neuronalwiring.html}{}.
\label{C-Elegans}

\subsection{\text{[Flagser]}} \url{https://github.com/luetge/flagser}{}.
\label{Flagser}

\subsection{\text{[Flag-Adapt]}} \url{https://github.com/JasonPSmith/flagser-adaptions}{}.
\label{Flag-Adapt}

\subsection{\text{[Neurotop]}} \url{http://neurotop.gforge.inria.fr/}{}.
\label{Neurotop}

\subsection{\text{[KONECT]}} \url{http://konect.uni-koblenz.de/}{}.
\label{KONECT}

\subsection{\text{[UPHESS]}} \url{https://hessbellwald-lab.epfl.ch/hessbellwald/}{}.
\label{UPHESS}

\bibliographystyle{alpha}
\bibliography{flagser}

\end{document}